\documentclass[12pt]{article}
\usepackage{amssymb}
\newcommand{\henk}{Henstock--Kurzweil}
\newcommand{\R}{{\mathbb R}}
\begin{document}
\hspace{-2cm}
\raisebox{12ex}[1ex]{\fbox{{\footnotesize
Preprint
February 4, 2008.  To appear in {\it American Mathematical Monthly}
}}}
\begin{center}
{\large\bf Review of \\{\it A garden of integrals},\\ by Frank E. Burk
(MAA, 2007)}\\
{\it Dedicated to the memory of Ralph Henstock (1923-2007).}
\vskip.25in
Erik Talvila\footnote{Supported by the
Natural Sciences and Engineering Research Council of Canada.
}\\ [2mm]
{\footnotesize
Department of Mathematics and Statistics \\
University College of the Fraser Valley\\
Abbotsford, BC Canada V2S 7M8\\
Erik.Talvila@ucfv.ca}
\end{center}
{\bf Executive summary:}  Well-written set of extended exercises
that take the reader through Cauchy, Riemann, Lebesgue, \henk, Wiener
and Feynman integration.  Requires a thorough grounding in
undergraduate analysis.   Anyone interested in the
chapters on Cauchy and Riemann integration 
probably won't have the prerequisites to read them.\\

\noindent
Riemann, Lebesgue, Denjoy, \henk,  McShane, Feynman, Boch-ner.  There
are well over 100 named integrals.  Why so many?  Some are
of historical interest 
and have been superseded by better, newer ones.  The Harnack
integral is subsumed by the Denjoy.  Some are equivalent, as are McShane
and Lebesgue in $\R^n$, and Denjoy, Perron, \henk\ in $\R$. 
Some are designed to work in special spaces:
Feynman for path integrals, where the domain is the set of all continuous
functions on $[0,1]$, Bochner for  Banach space-valued functions.  Others
are designed to invert special derivatives such as the symmetric derivative
or distributional derivative.  And, we continue
to keep the Riemann integral around because it is so easy to define.
Can't we have just one super-integral that does everything?  The
problem is that if an integral is defined to work, for example, 
on all Banach spaces
then it will probably be unnecessarily complicated when restricted to
a simple setting such as the real line.  There is also a pedagogical
issue.  We need to build up mathematics from simple bits  before we
can define the most abstract structures (Bourbaki notwithstanding).  
So, it seems we will have
to live with this plethora of integrals.

First there were Newton and Leibniz for whom integration was
anti-differentiation.  In the 1800s there were Cauchy and Riemann
who defined integrals in terms of approximating sums  formed
by partitioning the domain of a function on the real line.
There are two main themes to integration in the 20th century.
The first is the Lebesgue theory.  The Lebesgue integral of 1905
immediately revolutionised analysis.  Although more complicated
than the Riemann integral, it brought with it Borel's theory of 
measure.  The two concepts of measure and integral led to a very
robust theory of integration with powerful limit theorems (such
as Lebesgue dominated convergence), the ability
to extend the integral to abstract settings, and a Banach space of integrable
functions.  This allowed the flowering of such fields as probability and
potential theory.  

In the first decades of the 1900s, while the Lebesgue
theory was being refined by mathematicians such as Fatou and Radon,
there was a parallel development by Denjoy, Lusin and Perron.  The
Lebesgue integral is absolute:  function $f$ is integrable if and only
if $|f|$ is integrable.  But, functions such as $g(x) =x^2\sin(x^{-3})$
with $g(0)=0$ have a derivative that exists everywhere and for which
$|g'|$ is not integrable in any neighbourhood of the origin.  Thus
the Lebesgue integral cannot integrate all derivatives.  Integrals
introduced by Denjoy (1912) and Perron (1914) were able to integrate 
all derivatives.  But, these had rather cumbersome definitions, and
it was not easy to define them in settings other than the real line.
By World War II work on Denjoy--Perron
integration had more or less ceased.  
In the late 1950s Henstock and Kurzweil independently
defined an integral using a modified type of Riemann sum.  This integral
was equivalent to that of Denjoy and Perron and was much easier to
work with.  In fact, it's probably more transparent than the
improper Riemann integral.  It has reasonable limit theorems and  integrates all
derivatives.  Working with measures is also easy. 
Whereas the Lebesgue integral has difficulty dealing with integrals such
as $\int_0^\infty\sin(x^2)\,dx$, which exists as a conditionally
convergent improper Riemann integral, these are not problems
in the \henk\ theory. It was  for these kinds of reasons that the
\henk\ integral won many adherents and there are quite a number
of textbooks that discuss it at the undergraduate or graduate level,
such as \cite{bartle} and \cite{gordon}. But the \henk\
integral
lacks a Banach space structure and there is no canonical extension to
$\R^n$.  So its development 
won't be the end of the search for the perfect integral.

With so many integrals it is not surprising that many surveys of
integration have been published.  First there are the historical.
An excellent short summary of integration theory
is given by Henstock in the otherwise unreadable \cite{henstock}.
Methods due to Archimedes and his predecessors, and the integrals
of Cauchy and 
Riemann are discussed in Edwards' {\it The historical development
of the calculus} \cite{edwards}.
See Hawkins \cite{hawkins} for history of the Lebesgue integral.
For a history of non-absolute integration, see Bullen's survey \cite{bullen}.
The handbook by Zwillinger \cite{zwillinger} 
gives a brief description, an example, and
a few references for most types of integrals as well as many
techniques of integration, both analytical and numerical.  It's
an excellent starting point for anyone interested in integration.

Here are
five books that try to prove the major theorems for several 
different integrals.  First is {\it Modern theories of integration}
\cite{kestelman}
by Kestelman, who was one of Henstock's Ph.D. examiners.  Modern
means 1937. Working at perhaps a beginning graduate level, the
author gives a rather complete description of the Riemann integral, improper
Riemann integral, Lebesgue measure and Lebesgue integral in $\R^n$.
The Riemann--Stieltjes and Denjoy integrals are defined in $\R$.  
Along the same lines but including abstract measure and integration
as well as full descriptions of the Perron and Denjoy integrals is
the classic and authoritative {\it Theory of the integral} by Saks
\cite{saks}.  Pesin \cite{pesin} includes various specialised integrals.  Using Lebesgue measure in $\R^n$, {\it
Theories of integration} by Kurtz and Swartz \cite{kurtz} is
an undergraduate text that develops Riemann, Lebesgue and McShane
integration.  It assumes some background in analysis.
Thus, the reader is expected to know what {\it countable} means
but concepts such as {\it limsup} are defined.  It has lots of exercises
and is a well-crafted textbook.  At a slightly more advanced level
is Gordon's {\it The integrals of Lebesgue, Denjoy, Perron, and Henstock}.
It covers the integrals in the title as well as such topics as Darboux
and Baire class one functions.  It has good exercises and is a very carefully
written, well-crafted text at the graduate level.  One shortcoming
is that it only discusses integrals on $[0,1]$.  Nonetheless, it
has become a standard reference, displacing Saks.

With this background behind us,
let's get down to the business of {\it Garden of integrals}.
This book was written by Frank E. Burk, who is author of the successful
text {\it Lebesgue measure and integration: an introduction} \cite{burk}.
The author is now deceased. A pity as the book should have good legs.
Not everything I'm going to say about the book is positive and I feel
a bit bad that the author won't have the chance to respond publicly.

I should tell you where I'm coming from.  I am familiar with Lebesgue
and \henk\ integration at the level of having written several research
papers on these integrals.  I'll be making a shameless plug for these
later.  I've taught calculus plenty of times.  I was almost completely
ignorant of the Wiener and Feynman integrals.  So, some of the chapters
I was figuring would be old hat and from others
I was hoping  to learn some new stuff.

This book surveys several types of integrals.  
There are chapters on the Cauchy, Riemann, Riemann--Stieltjes,
Lebesgue, Lebesgue--Stieltjes,
Henstock--Kurzweil, Wiener and Feynman
integrals.  There is also an historical overview and a chapter
on Lebesgue measure.  The emphasis is purely theoretical with
very little in the way of applications.  Throughout the book the author assumes
a thorough grounding in undergraduate real analysis.  Thus
the book is pitched at senior undergraduate or beginning graduate
level.

The book is essentially a string of exercises, implicit or explicit,
in which the author has you work out the details of the various
theories.  So, this is definitely not a coffee table book.  To
read it you'll have to turn the TV off,  sit down, shut up and
get to work.  Working through it will teach you a lot of integration.

The writing style is attractive. Sentences flow well, and things
are laid out
logically.  In this sense it's a very good book.  What's not good
about it is that I'm not sure who it was written for.  As any
writing instructor will tell you, when you sit down to write 
the first item to have straight is: Who is your audience?  I doubt
that was done in this case.  There are beginning chapters on 
Cauchy and Riemann integrals.  The Cauchy integral is like the
Riemann except that in your Riemann sum  you always use the left endpoint 
of each subinterval when evaluating the integrand.
These are pretty basic integrals but to read these chapters
you need to know the 
following: trapezoid rule, integration by parts, improper
Riemann integral, Cauchy sequences,
uniform convergence, Weierstrass $M$-test,
Fourier series, separation of variables for the Laplace equation,
Cantor set, nowhere dense.  That's a tall order.  I'm wondering
just who it is who would know all this stuff but still need to
learn something about the Cauchy and Riemann integrals.  As I was
reading these chapters I was very puzzled as to who the book was for.
Then on page 65 (Riemann integral) 
the author lets the cat out of the bag, ``I dare say
we can all remember our first encounter with [the Cantor set]."
So, to read the book you need to have many results of 
real analysis at your fingertips.  
That said, these chapters have some interesting
material and some good examples,  such as the construction of
a continuous nowhere differentiable function and
a construction of
a Riemann integrable function with a dense set of discontinuities.

As Cauchy did, the author defines the Cauchy integral only
for continuous functions.  Burk is correct in proclaiming that in this sense 
the Riemann integral
properly includes the Cauchy integral.  However, as already
pointed out by my editor, 
the space of all Cauchy integrable
functions (dropping the continuity assumption) 
is the same as the space of Riemann integrable
functions.  This was proved by Gillespie \cite{gillespie}.
Burk also states that indefinite Cauchy and Riemann integrals
are absolutely continuous, but absolute continuity isn't defined
for another 65 pages.  It would have been better to give the
easy proof that these integrals are Lipshitz continuous.

Preceding these chapters is an historical survey.  It has
some discussion on geometric methods of the ancients
and
then brief descriptions of the integrals covered in later
chapters.  These were too vague to be of much interest and this
chapter should probably have been omitted.
A better list of historical sources would have been appreciated.

The chapters on Lebesgue measure and Lebesgue integration are much 
better. Voltera's example of a function with a bounded derivative that
is not Riemann integrable makes good motivation for introducing
the Lebesgue integral.
These chapters bring in many of the standard topics in
measure and integration on the real line.
Most concepts here are followed by useful exercises
although the coverage of sigma algebras and Borel sets is very brief.
Bounded variation and absolute continuity are defined, but there
is so little discussion of them I'm not sure if the uninitiated would
get the whole picture.  Again, you need to know your real analysis
to make your way through this book.  Thorough familiarity with topics
such as {\it inf}, {\it sup} and the Heine--Borel theorem are assumed
here.

The short chapters
on  Riemann--Stieltjes
and Lebesgue--Stieltjes integrals don't add much and could well have
been omitted.

One of the themes of the book is the Fundamental Theorem of Calculus.
At the end of each chapter is a summary that consists of the Fundamental
Theorem as it applies to that integral.  This is all motivation for
the \henk\ integral, since it has the best version of this theorem on
the real line.  With the Lebesgue integral we know that the primitives
are the absolutely continuous functions.  The corresponding class of
functions for the \henk\ integral is more complicated.  It's called
{\it generalized 
absolute continuity in the restricted sense} ($ACG*$,
see \cite{gordon}).  A reasonable approach that some authors use is
to take the continuous functions that are differentiable except for
countable sets as the primitives.  This is done here,  but no
mention of the real set of primitives is made.  A problem with
the \henk\ integral is that there is no canonical way to extend 
$ACG*$ to $\R^n$, let alone to more abstract spaces,  and the
set of integrable functions forms an incomplete normed space.
We can get around both of these problems by using $C^0$
as the space of primitives and saying the integral inverts the
distributional derivative, so that the integrable distributions
are those that are the distributional derivative of continuous
functions.  The space of integrable distributions
is then  a Banach space that includes $L^1$ and the \henk\ integrable
functions.
This is a really simple way to define a very general integral.
Using the divergence theorem, this definition can be used in $\R^n$.
See \cite{talvila} for an easy but complete description of this integral.

The chapter on \henk\ integration is good.  Initially it proceeds at
a much more gentle pace than the Lebesgue chapters.  The author takes
the time to motivate this integral and then works out its properties
in some detail.

The final two chapters are the tour de force.  In treating the
Wiener and Feynman integrals in a book of this level (and length)
he's taking on a big task and breaking new ground.  
The domain of integration is now infinite dimensional:  
the set of continuous real-valued functions on $[0,1]$.
These chapters require Lebesgue measure and integration as well
as some familiarity with probability, the heat equation, complex
variables, and linear operators on a Hilbert space.  They're written
at a higher level than the preceding chapters.  The Wiener integral
is described reasonably well. The examples and exercises should
help you understand what's going on.

I was especially impressed with
the Feynman integral chapter.  It is long and difficult
and packed with  lots of stuff: Schr\"odinger equation,
Fourier transform in the Schwartz space, Trotter product formula, semigroups of
linear operators.  Many of these topics are dealt with quite briefly,
so working out the details and filling the gaps will take time and
consultation of outside references.   I can't say I understood it all
on a first reading but that gives me plenty to go back to.  

I found it
irksome that Plank's constant is never mentioned by name and is defined  
as $\hbar=1.054 \times 10^{-27}$$\,erg\,sec$.  
I never knew it was a rational number.

A note regarding peripherals.  The typesetting is generally good and the book is
easy to read, with lots of white
space and well-displayed formulas.
Each chapter begins with an irrelevant quotation by some famous
person.  I suppose these are meant to inspire.  They didn't do anything for me.
For witty and inspiring beginning of chapter quotes, I much prefer
\cite{reed}.
There are about 50 illustrations,
mostly in the historical overview, Wiener and Feynman chapters.
These are quite good.
The index is thin.
There aren't many typos, although one
recurring one was quotes written "hey there" and not ``hey there".
The editor should have caught this.  The binding is good quality,
bound in signatures.  At \$51.95 (cheaper for MAA members) the
book is excellent value.  The cover photo is a bonus.

A foreword of a book is usually written by someone other than the author.
It is generally some sort of endorsement.
Publishers often have someone of authority write it and then set their
name in larger type than the author.  I once saw a book with a preface
by Eric Idle of Monty Python fame.  He basically just wrote: Here's the
foreword now send me my cheque.  If you're an unknown environmentalist
and they get Al Gore to write the foreword to your book then you're in
good shape.  Meanwhile, a preface is usually written by the author and it will
tell you something about the book.  Some authors have written prefaces
that have become famous in their own right.  In the second edition
of {\it Lord of the Rings} (1966), 
J.R.R. Tolkien makes it clear that the book was
not intended as an allegory of World War II.  For {\it The art of
computer programming}, Donald Knuth includes a flowchart that has
become legendary \cite{knuth}.

Mathematics books usually include a preface that tells the reader
roughly what topics will be covered, the chapter dependencies,
the level of depth, and the
required prerequisites.  {\it Garden of integrals} has no preface.  
It has a short and
completely uninformative
foreword written by the author.  This is unfortunate as some 
unsuspecting souls may try to read this attractive looking book
and find it pretty tough going.  It is a major
failing of the author/editor to not have told the prospective reader
what they are expected to know in order to make sense of any given
chapter of this book.   Again, for every chapter the 
prerequisites are a thorough knowledge of real analysis
at the undergraduate level.

Topics I would have liked to have seen are first and foremost,
some discussion of the integral on sets other than compact intervals
on the real line.
Improper Riemann integrals are a glaring omission.
The discussion of the Riemann integral confines itself to $[a,b]$
but in the chapter on the Wiener integral, improper Riemann
integrals are assumed.  Integration in $\R^n$, Stokes theorem and 
some applications to probability would have been nice.

Another thing that bothered me were the references.  There is a list
at the end of each chapter.  These may have been the works (mostly
textbooks) the
author used to write the book but a more comprehensive list suitable
for further reading 
should have been compiled.

How does this book compare with the surveys mentioned above?
Except for the chapters on Wiener and Feynman integration the
book is comparable to Kurtz and Swartz \cite{kurtz} and Gordon
\cite{gordon}, although these
authors do a better job of incorporating material from the real 
analysis course.  All of the books cited above are pretty formal.
Burk strikes a nice balance between the formal definition, theorem,
proof style and chatty, informal discussion. He's more fun and
free wheeling.  He's engaging because he's constantly forcing
you to fill in little details.  This is a lousy book to use as
a reference to look up something like
the Lebesgue dominated convergence theorem,
but if you work through the book you'll have a good time learning
the proof.  And, his chapters on Wiener and Feynman integration
set him apart from the rest.

Who can benefit from this book? I would be happy to give it to good students
with the necessary background for use in an undergraduate reading course or
seminar.  Such a course could be on Lebesgue or \henk\ integration.
Since the book is pretty much just a long set of exercises,
this could work well.  However, many topics are dealt with in
a cursory manner so outside reading would have to be prescribed.
A graduate seminar could be set up
for the Wiener and Feynman chapters.  And, of course, the book is
ideal for self study.
Will I recommend my college library purchase this book?  Definitely.

\end{document}